\documentclass[letterpaper, 10 pt, conference]{ieeeconf}  
\IEEEoverridecommandlockouts    
\usepackage{graphicx}

\title{\LARGE \bf
Nearly Optimal Chaotic Desynchronization of Neural Oscillators
}

\author{Jeff Moehlis, Michael Zimet, and Faranak Rajabi% <-this % stops a space
\thanks{J. Moehlis is with the Department of Mechanical Engineering and the Program in Dynamical Neuroscience, 
        University of California, Santa Barbara, United States of America
        {\tt\small moehlis@ucsb.edu}}%
\thanks{M. Zimet is with the Program in Dynamical Neuroscience, 
        University of California, Santa Barbara, United States of America
        {\tt\small michaelzimet@ucsb.edu}}%
\thanks{F. Rajabi is with the Department of Mechanical Engineering, 
        University of California, Santa Barbara, United States of America
        {\tt\small faranakrajabi@ucsb.edu}}%        
}

\begin{document}

\maketitle
\thispagestyle{empty}
\pagestyle{empty}

%%%%%%%%%%%%%%%%%%%%%%%%%%%%%%%%%%%%%%%%%%%%%%%%%%%%%%%%%%%%%%%%%%%%%%%%%%%%%%%%
\begin{abstract}
Motivated by deep brain stimulation treatment of neural disorders such as Parkinson's disease, it has been proposed that desynchronization of neural oscillators can be achieved by maximizing the Lyapunov exponent of the phase difference between pairs of oscillators.  Here 
we consider two approximations to optimal stimuli for chaotic desynchronization of neural oscillators.  These approximations are based on the oscillators' phase response curve, and unlike previous approaches do not require numerical solution of a two-point boundary value problem.  It is shown that these approximations can achieve nearly optimal desynchronization, and can be used with an event-based control scheme to desynchronize populations of noisy, coupled neurons.

\end{abstract}

%%%%%%%%%%%%%%%%%%%%%%%%%%%%%%%%%%%%%%%%%%%%%%%%%%%%%%%%%%%%%%%%%%%%%%%%%%%%%%%%
\section{INTRODUCTION}

Evidence suggests that pathological synchronization of neurons in the motor control region of the brain is a factor which contributes to symptoms such as tremors for patients with Parkinson's disease~\cite{hamm07}.  This has led to the hypothesis that Deep Brain Stimulation (DBS), which has been shown to be effective for treating the symptoms of Parkinson's disease, might be desynchronizing the neural activity~\cite{wils22},~\cite{tass07}, \cite{wils11}.  

The standard implementation of DBS uses open loop injection of a pulsatile, high frequency signal.  There have been various attempts to design alternative control methods to achieve this goal using both closed-loop and open-loop approaches; see~\cite{wils22} for a review. Many of these approaches attempt to minimize the amount of energy which is necessary to desynchronize the neural population; this could extend the battery life of the pacemaker which drives DBS, and could also minimize side effects from the procedure because weaker inputs are used.
For example,~\cite{wils14a} proposed the use of chaotic desynchronization, in which the control is chosen to maximize the Lyapunov exponent associated with phase differences for the neurons.  An alternative approach described in~\cite{danz10,nabi12} proposes the use of optimal phase resetting, in which an input drives the system to the phaseless set where the neurons are particularly sensitive to noise.  This approach has recently been extended to a stochastic control framework in~\cite{raja25}. References~\cite{mong18acc,mong19} used phase density control in which an input drives the system toward a desired phase distribution, which can be a flat distribution corresponding to complete desynchronization, or a distribution corresponding to clustering, in which each cluster of neurons fires out of phase with other clusters.  Other notable approaches to desynchronization and clustering include 
\cite{rose04}-\nocite{popo06}\nocite{rose04pre}\nocite{scho09}\nocite{kiss07}\nocite{Lucken2013}\nocite{Lysyansky2011}\nocite{holt16a}\nocite{wils15cluster}\nocite{matc18}\nocite{wils16}\nocite{matc21}\nocite{wils20}\nocite{qin23}\cite{vu24}.

In this work, we revisit the optimal chaotic desynchronization method from~\cite{wils14a}.  That paper calculated the optimal input by numerically solving the two-point boundary value problem which arises from the optimal control problem.  Here we instead explore approximations that can be found directly from the neuron's phase response curve, which characterizes how the timing of a stimulus affects the neuron's dynamics.  We show that one approximation in particular can give nearly optimal desynchronization.  Moreover, we show that the approximations can be implemented with an event-based control scheme to desynchronize populations of noisy, coupled neurons.  We emphasize that the approximations are found with less computational effort than the approach in~\cite{wils14a}, and also provide some insight into why the optimal input has the form which it does.  We also show that the approximations can be used in the related problem of synchronization.

The organization of the paper is as follows.  In Section \ref{lyapsection}, we derive expressions for two approximations to the input which achieves optimal chaotic desynchronization. Section \ref{desync} investigates how well these inputs work for desynchronizing pairs of oscillators for three different phase response curves.  Section~\ref{pop} demonstrates the use of the approximately optimal stimuli for desynchronizing populations of noisy, coupled neurons.  Section \ref{sync} considers the related problem of optimal synchronization.  Concluding remarks are given in Section \ref{conclusion}.

\section{Derivation of Optimal Control Input and its Approximations} \label{lyapsection}
We first review the procedure from~\cite{wils14a} for finding an energy-optimal stimulus which exponentially desynchronizes a population of neural oscillators (that is, neurons which fire periodically in time in the absence of control, noise, and coupling).  This approach is based on the neurons' phase response curve (PRC), a function which captures the response properties of an individual neuron to a stimulus.  This can in principle be measured experimentally~\cite{neto12}, or by solving the appropriate adjoint equation if the equations and parameters in the neural model are known, e.g.~\cite{tutorial}. Phase reduction of the dynamics give the following model for a single neuron~\cite{brow04,tutorial}: 

\begin{equation}\label{PRCeq}
\frac{d\theta}{dt}=\omega+Z(\theta)u(t),
\end{equation}
where $\theta \in [0,2 \pi)$ is the phase of the neuron, with $\theta=0$ corresponding to the spiking of the neuron in our convention, and the PRC $Z(\theta)$ is $2 \pi$-periodic.  The neuron's angular frequency is $\omega = 2 \pi / T$, where $T$ is the neuron's natural period. The control input $u(t) = I(t)/C$, where $I(t)$ is the input current and $C=1\mu\mathrm{F/cm}^2$ is the constant neural membrane capacitance.  

Now consider two identical neurons (i.e., with the same $\omega$ and $Z(\cdot)$) which are subjected to a common stimulus $u(t)$:
\begin{equation} \label{basedyn}
\frac{d\theta_i}{dt} = \omega + Z(\theta_i)u(t), \quad i = 1,2.
\end{equation}
Assuming that the neurons are nearly in-phase, so that $\theta_1 \approx \theta_2$, and letting $\phi \equiv  \theta_2-\theta_1 $, we obtain
\begin{equation}\label{phieq}
\frac{d\phi}{dt}=Z'(\theta)u(t)\phi + \mathcal{O}(\phi ^2),
\end{equation}
where we have replaced $\theta_1$ by $\theta$.
Neglecting $\mathcal{O}(\phi^2)$ and higher order terms in (\ref{phieq}),
%\begin{equation}
 $   \phi(t) = \phi_0 e^{\int_0^t Z'(\theta(s)) u(s) ds}.$
%\end{equation}
Thus, $\phi(T) = \phi_0 e^{\Lambda T}$, where $\Lambda$ is the finite time Lyapunov exponent, here defined over one natural period of the oscillator (cf. \cite{abou09}):
\begin{equation} \label{lyapexp2}
%\Lambda = \frac{\log(\phi(T))}{T}  = \frac{\int_0^T Z'(\theta(s))u(s)ds}{T}.
\Lambda \equiv \frac{\int_0^T Z'(\theta(s))u(s)ds}{T}.
\end{equation}

%\section{Non-Charge-Balanced (NCB) Control}

Now consider the following optimal control problem.  Suppose for all non-zero\footnote{If $u=0$ for all $0\le t \le T$, then $\Lambda=0$, thus no desynchronization.} stimuli $u(t)$ that steer $\theta$ from $\theta(0)=0$ to $\theta(T)=2 \pi$ (that is, stimuli that do not create any net change of phase over the natural period of the neuron), we want to find the stimulus that minimizes the cost function 
\begin{equation}
G[u(t)]=\int_0^{T} [u(t)^2 - \beta Z'(\theta(t))u(t)]dt.
\label{G}
\end{equation}
%\footnote{If $u=0$ for all $0\le t \le T$, then $\Lambda=0$, thus no desynchronization.}
Here, $\int_0^{T} [u(t)^2]dt$ is the energy associated with the stimulus, and $\beta>0$ is a weighting parameter that determines the relative importance of maximizing the Lyapunov exponent versus minimizing the energy.  Thinking of this as an optimal control problem, we seek to minimize~\cite{forg04}
\begin{eqnarray}
    \label{costfunction}
\mathcal{C}[u(t)] &=& \int_0^{T}  \left\{    u(t)^2 - \beta Z'(\theta)u(t) \right. \nonumber \\
&& \left. + \lambda \left(    \frac{d\theta}{dt} - \omega -Z(\theta)u(t)  \right)   \right\} dt.
\end{eqnarray}
Here $\lambda$ is a Lagrange multiplier which enforces the neural dynamics obey~(\ref{PRCeq}).  The corresponding Euler-Lagrange equations are~\cite{wils14a}

\begin{eqnarray}
\dot{\theta}\!\!\!\!&=&\!\!\!\!  Z(\theta) \left[ \beta Z'(\theta) +\lambda Z(\theta)\right]/2  +\omega, \label{thetadot}\\
\dot{\lambda} \!\!\!\!&=& \!\!\!\! - \left[ \beta Z'(\theta) +\lambda Z(\theta)\right]  \left[ \beta Z''(\theta)+\lambda Z'(\theta) \right]/2, \label{lambdadot} \\
u(t) \!\!\!\! &=&\!\!\!\! [\beta Z'(\theta) +\lambda Z(\theta)]/2 ,\label{findcontrol}
\end{eqnarray}
where $\dot{(~)} = d(~)/dt$ and $(~)'=d(~)/d\theta$.  Reference~\cite{wils14a} uses a shooting method to numerically solve (\ref{thetadot}) and (\ref{lambdadot}) subject to the boundary conditions $\theta(0)=0$, $\theta(T)=2 \pi$.  From this, one can then find the optimal input $u^*(t)$ from (\ref{findcontrol}).
 
In this paper, we instead consider approximate solutions to (\ref{thetadot}-\ref{findcontrol}).  We observe that because we seek to minimize the input energy $\int u(t)^2 dt$, the control input $u(t)$ will be relatively small.  From (\ref{PRCeq}), in the limit of vanishing $u(t)$,
\begin{equation}
    \theta(t) \approx \omega t.
\label{theta_approx}
\end{equation}
Next, for $\beta$ sufficiently large we assume that for solutions to the two-point boundary value problem (\ref{thetadot}-\ref{findcontrol}), $|\lambda Z(\theta)|$ is small compared with $|\beta Z'(\theta)|$; this is supported by the numerical results in~\cite{wils14a}.
Using this and (\ref{theta_approx}) in~(\ref{findcontrol}), we obtain the following approximation for the optimal input:
\begin{equation}
    u(t) \approx \frac{\beta}{2} Z'(\omega t) \equiv \tilde{u}_1^*(t).
\end{equation}

%\begin{eqnarray}
%    \dot{\theta} &\approx& \frac{\beta}{2} Z(\theta) Z'(\theta) + \omega \\
%    &\approx& \frac{\beta}{2} Z(\omega t) Z'(\omega t) + \omega \\
%    &=& \frac{\beta}{4 \omega} \frac{d}{dt} \left[ Z(\omega t) \right]^2 + \omega.
%    \end{eqnarray}
%\end{eqnarray}
We next assume that $|\lambda Z'(\theta)|$ is small compared with $|\beta Z''(\theta)|$.  Under this and the above assumption,
%\begin{eqnarray*}
%    \dot{\lambda} &\approx& -\frac{\beta^2}{2} Z'(\theta) Z''(\theta) \\
%    & \approx & -\frac{\beta^2}{2} Z'(\omega t) Z''(\omega t) \\
%    &=& -\frac{\beta^2}{4 \omega} \frac{d}{dt} \left[ Z'(\omega t) \right]^2.
%\end{eqnarray*}
\begin{eqnarray*}
    \dot{\lambda} &\approx& -\frac{\beta^2}{2} Z'(\theta) Z''(\theta) \approx  -\frac{\beta^2}{2} Z'(\omega t) Z''(\omega t) \\
    &=& -\frac{\beta^2}{4 \omega} \frac{d}{dt} \left[ Z'(\omega t) \right]^2.
\end{eqnarray*}
\begin{equation}
  \Rightarrow \qquad  \lambda(t) \approx -\frac{\beta^2}{4 \omega} \left[ Z'(\omega t) \right]^2 + c,
    \label{lambda_approx1}
\end{equation}
where $c$ is a constant.
We now argue that to this order of approximation, the Lyapunov exponent does not depend on $c$.  In particular, from (\ref{findcontrol})
\begin{eqnarray*}
u(t) &=& \frac{1}{2} \left[ \beta Z'(\theta) + \lambda Z(\theta) \right] \approx  \frac{\beta}{2} Z'(\omega t) + \frac{1}{2} \lambda Z(\omega t) \\
&\approx& \frac{\beta}{2} Z'(\omega t) - \frac{\beta^2}{8 \omega} \left[ Z'(\omega t) \right]^2 Z(\omega t) + \frac{c}{2} Z(\omega t).
\end{eqnarray*}
Plugging this into (\ref{lyapexp2}), we find the following approximate expression for the Lyapunov exponent:
\begin{eqnarray}
    \Lambda \!\! &\approx& \!\! \frac{1}{T} \int_0^T Z'(\omega s) \left[\frac{\beta}{2} Z'(\omega s) - \frac{\beta^2}{8 \omega} \left[Z'(\omega s) \right]^2 Z(\omega s) \right] ds \nonumber \\
    && + \frac{c}{2 T} \int_0^T Z'(\omega s) Z(\omega s) ds. \label{lyap_approx}
\end{eqnarray}
But the last integral in (\ref{lyap_approx}) vanishes:
\begin{eqnarray*}
    \int_0^T Z'(\omega s) Z(\omega s) ds &=& \frac{1}{2 \omega} \int_0^T \frac{d}{ds} \left[ Z(\omega s) \right]^2 ds \\
    &=& [Z(2 \pi)]^2 - [Z(0)]^2 = 0,
\end{eqnarray*}
where we have used the fact that $\omega T = 2 \pi$, and that the PRC is $2 \pi$-period.  Thus, we can choose $c$ to be any convenient value.  For simplicity in the following, we choose $c = 0$.  Then (\ref{lambda_approx1}) simplifies to 
\begin{equation}
    \lambda(t) \approx -\frac{\beta^2}{4 \omega} \left[ Z'(\omega t) \right]^2.
    \label{lambda_approx2}
\end{equation}
Using this in (\ref{findcontrol}), we obtain the following approximation to the optimal input:
\begin{equation}
    u(t) \approx \frac{\beta}{2} Z'(\omega t) - \frac{\beta^2}{8 \omega} \left[ Z'(\omega t) \right]^2 Z(\omega t) \equiv \tilde{u}_2^*(t).
    \end{equation}

We note that the energy associated with the stimulus $u^*(t)$ will in general not be the same as the energy associated with the stimuli $\tilde{u}^*_1(t)$ or $\tilde{u}^*_2(t)$.  In order to make fair comparisons between these inputs, in the following we rescale $\tilde{u}^*_1(t)$ and $\tilde{u}^*_2(t)$ so that all inputs use the same energy for a given value of $\beta$.
In particular,
\[
u^*_i(t) = \tilde{u}^*_i(t) \times \sqrt{\frac{\int_0^T [u^*(t)]^2 dt}{\int_0^T [\tilde{u}^*_i(t)]^2 dt}}, \qquad i = 1,2.
\]
With this re-scaling,
\[
\int_0^T [u^*(t)]^2 dt = \int_0^T [u^*_1(t)]^2 dt = \int_0^T [u^*_2(t)]^2 dt.
\]

\section{Desynchronization Results} \label{desync}
\begin{figure}
    \centering
    \includegraphics[width=2.9in]{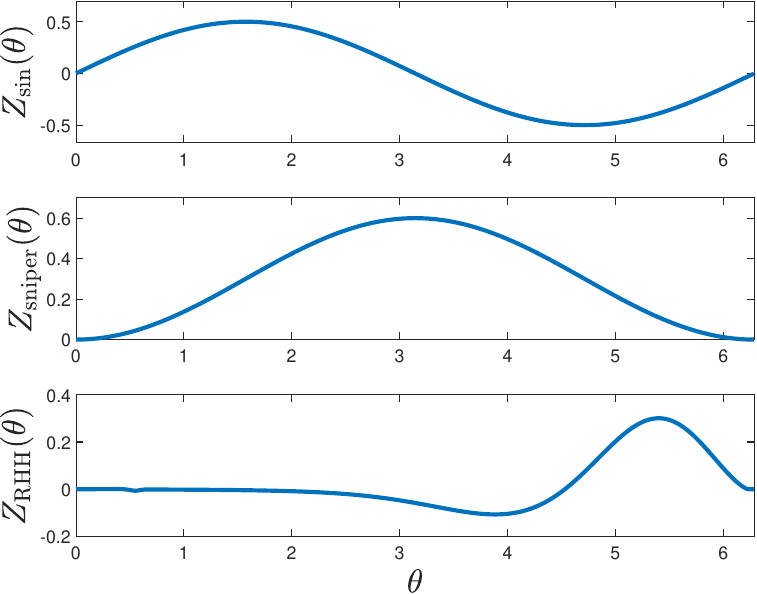}
    \caption{PRCs $Z_{\sin}(\theta)$, $Z_{\rm sniper}(\theta)$, and $Z_{\rm RHH}(\theta)$.}
    \label{fig:prcs}
    \end{figure}

We investigate the performance of the approximations $u^*_1(t)$ and $u^*_2(t)$ for three different examples PRCs: (i) $Z_{\sin}(\theta) \equiv 0.5 \sin(\theta)$, (ii) $Z_{\rm sniper}(\theta) \equiv 0.3 (1 - \cos \theta)$, and (iii) $Z_{\rm RHH} (\theta)$ for the Reduced Hodgkin-Huxley (RHH) equations~\cite{keen98} given by
\begin{equation}
\dot{V} = f_V(V, n), \qquad \dot{n} = f_n(V, n).
\label{eq:Vndot}
\end{equation}
Here $V$ is the voltage across the membrane, and $n$ is a gating variable.  The functions $f_V$ and $f_n$ and parameter values can be found in the Appendix.  We use a baseline current of $10 \mu A/cm^2$, which gives a natural period $T = 11.85$ ms.
Note that the PRCs $Z_{\sin}$ and $Z_{\rm sniper}$ could arise from a dynamical system near a Hopf bifurcation or a SNIPER (also known as a SNIC) bifurcation, respectively~\cite{brow04}.  The PRC $Z_{\rm RHH}$ was found numerically by solving the appropriate adjoint equation with XPP~\cite{erme02}, and is approximated as a Fourier series with the first 200 $\sin(\cdot)$ and $\cos(\cdot)$ terms. The PRCs considered in this paper are shown in Figure~\ref{fig:prcs}.

\begin{figure}
    \centering
    \includegraphics[width=2.9in]{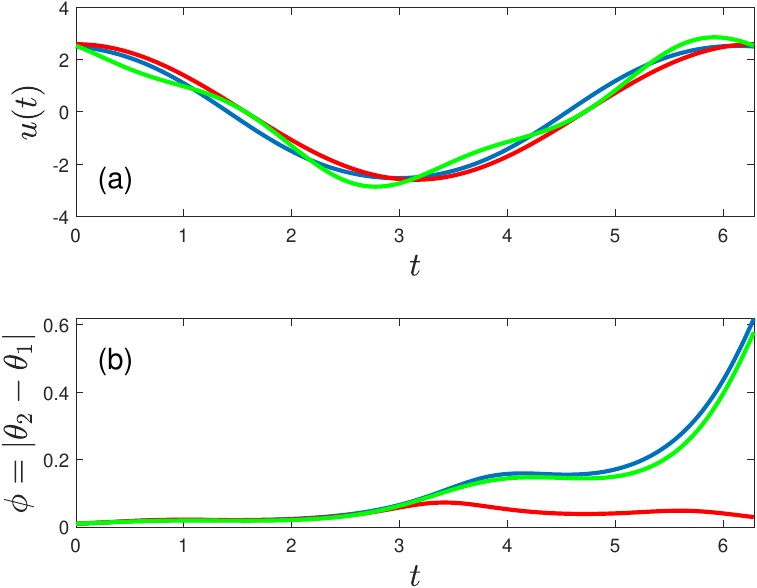}
    \caption{Results for the sinusoidal PRC $Z_{\sin}(\theta) = 0.5 \sin \theta$, with $\omega = 1$ and $\beta = 10$.  (a) shows the inputs $u^*(t)$ (blue), $u^*_1(t)$ (red), and $u^*_2(t)$ (green).  (b) shows the evolution of the phase difference $\phi$ for one cycle of each input, with colors as for (a).  The initial phase separation $\phi_0$ = 0.01.}
    \label{fig:sin}
    \end{figure}

\begin{figure}
    \centering
    \includegraphics[width=2.9in]{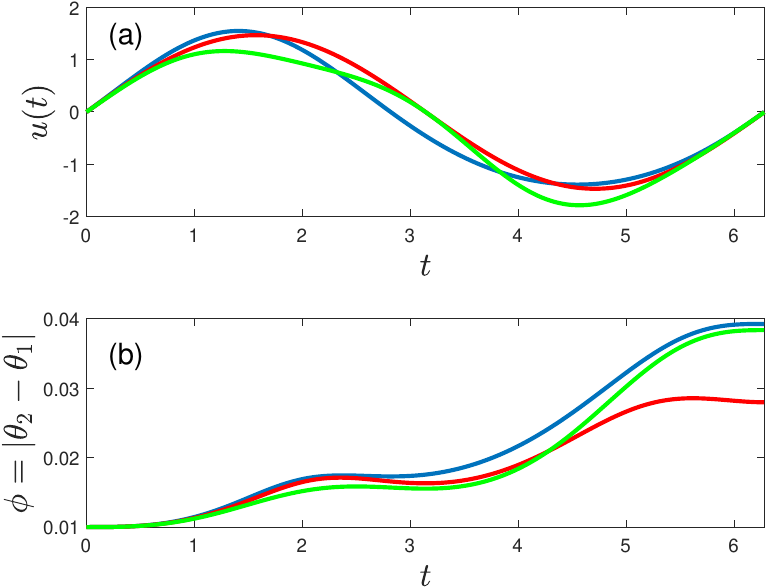}
    \caption{Results for the SNIPER PRC $Z_{\rm sniper}(\theta) = 0.3 (1 - \cos \theta)$, with $\omega = 1$ and $\beta = 10$.  (a) shows the inputs $u^*(t)$ (blue), $u^*_1(t)$ (red), and $u^*_2(t)$ (green).  (b) shows the evolution of the phase difference $\phi$ for one cycle of each input, with colors as for (a).  The initial phase separation $\phi_0$ = 0.01.}
    \label{fig:sniper}
    \end{figure}

 The computed inputs for $Z_{\sin}(\theta)$, $\omega = 1$, and $\beta = 10$ are shown in Figure~\ref{fig:sin}(a).  We note that the $L_2$ distance between $u^*(t)$ and $u^*_2(t)$ is smaller than the $L_2$ distance between $u^*(t)$ and $u^*_1(t)$.
 More importantly, by looking at the final phase differences which result from these inputs, as shown in (b) where $\theta_1(0) = 0$ and $\theta_2(0) = 0.01$, we see that $u^*_2(t)$ gives larger final phase separation than $u^*_1(t)$.  In fact, the performance of $u^*_2(t)$ is very close to the performance of $u^*(t)$; thus, $u^*_2(t)$ is nearly optimal.
Results from a similar analysis for $Z_{\rm sinper}(\theta)$, $\omega = 1$, and $\beta = 10$, shown in Figure~\ref{fig:sniper}, are qualitatively similar.  
Interestingly, the $L_2$ distance between $u^*(t)$ and $u^*_2(t)$ is now larger than the $L_2$ distance between $u^*(t)$ and $u^*_1(t)$.
However, the final phase separation which results from $u^*_2(t)$, as shown in Figure~\ref{fig:sniper}(b) where $\theta_1(0) = 0$ and $\theta_2(0) = 0.01$, is nearly the same as that which results from $u^*(t)$, whereas the final phase separation which results from $u^*_1(t)$ is lower.  Again, $u^*_2(t)$ is nearly optimal.

\begin{figure}
    \centering
    \includegraphics[width=2.9in]{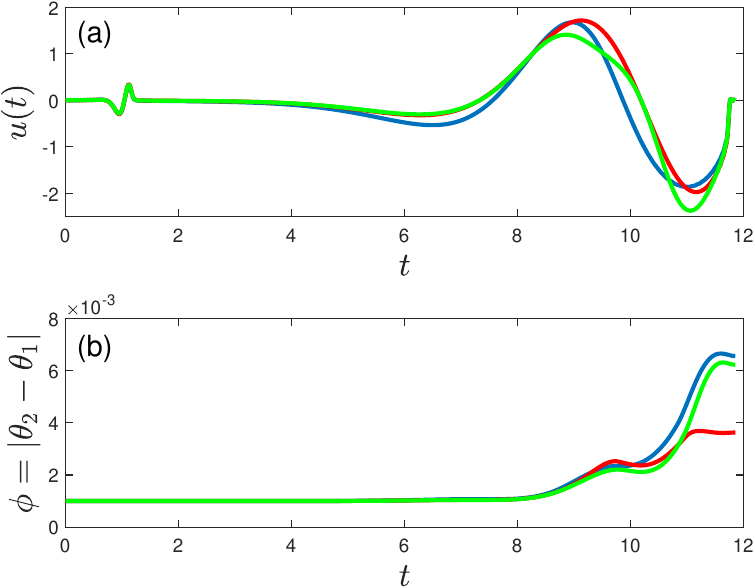}
    \caption{Results for the RHH PRC $Z_{\rm RHH}(\theta)$, with $\omega = 2 \pi/11.85$ and $\beta = 7$.  (a) shows the inputs $u^*(t)$ (blue), $u^*_1(t)$ (red), and $u^*_2(t)$ (green).  (b) shows the evolution of the phase difference $\phi$ for one cycle of each input, with colors as for (a).  The initial phase separation $\phi_0$ = 0.001.}
    \label{fig:RHH}
    \end{figure}

Figure~\ref{fig:RHH} shows results for the PRC $Z_{\rm RHH}(\theta)$ for the RHH equations, with (a) showing the inputs, and (b) showing the phase difference which results from these inputs.  We choose $\beta = 7$.  We find that the $L_2$ distance between $u^*(t)$ and $u^*_2(t)$ is larger than the $L_2$ distance between $u^*(t)$ and $u^*_1(t)$. 
However, again, $u^*_2(t)$ gives nearly optimal performance, and better performance than $u^*_1(t)$.

\begin{figure}
    \centering
    \includegraphics[width=2.9in]{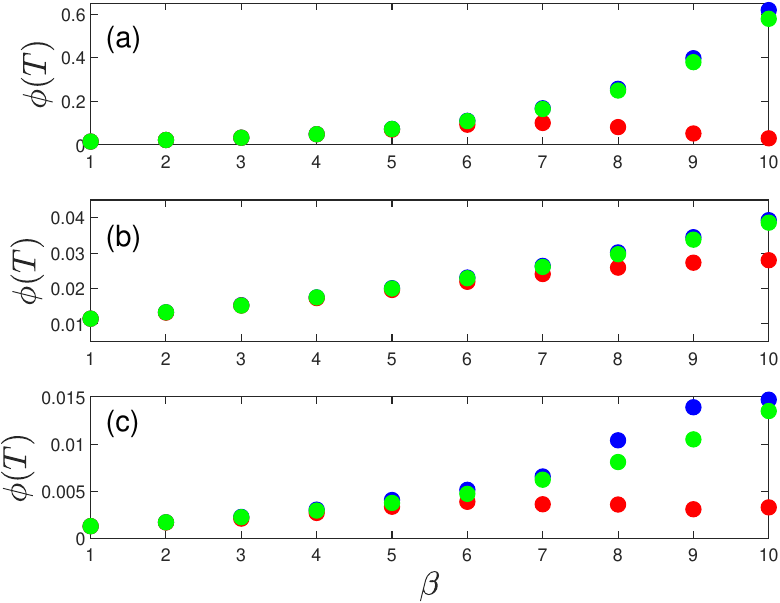}
    \caption{Final phase difference values $\phi(T)$ after one cycle of the different inputs, as a function of $\beta$.  The dots represent results for $u^*(t)$ (blue), $u^*_1(t)$ (red), and $u^*_2(t)$ (green).  The approximation $u^*_2(t)$ is nearly optimal for a range of $\beta$ values.}
    \label{fig:delta_theta}
    \end{figure}

In Figure~\ref{fig:delta_theta}, we show results for different values of $\beta$ for the three different PRCs.  Here we just show the final phase difference at the end of one cylce of input, that is, $\phi(T)$.  It is clear from this figure that the input $u^*_2(t)$ is nearly optimal across a range of $\beta$ values, whereas $u^*_1(t)$ is far from optimal for larger $\beta$ values.

\section{Application to Neural Populations} \label{pop}
While the results in the previous section illustrate that $u^*_2(t)$ gives nearly optimal performance when applied to two uncoupled oscillators, we are also interested in investigating the application of this stimulus to a larger population of coupled neural oscillators.  In this section, we consider a population of RHH neurons with all-to-all electrotonic coupling, and independent additive noise for each neuron.  The governing equations are:
\begin{eqnarray*}
\dot{V}_i &=& f_V(V_i, n_i) + \frac{1}{N}\sum_{j=1}^{N} \alpha (V_j - V_i) + u(t) + \eta_i(t),\\
\dot{n}_i &=& f_n(V_i, n_i),
\label{eq:Vndot_pop}
\end{eqnarray*}
where $N = 100$ is the number of neurons, $u(t)$ is the common input for all neurons, $\alpha = 0.04$ is the coupling strength, and \( \eta_i(t) = \sqrt{2D} \mathcal{N}(0,1) \) is intrinsic noise modeled as zero-mean Gaussian white noise with variance \( 2D \), where $D = 2$.  We emphasize that each neuron receives the same input $u(t)$, but a different realization of noise $\eta_i(t)$.  It is useful to think of the coupling as having a synchronizing effect, and the noise as having a desynchronizing effect.

To test our inputs, we use an event-based control scheme similar to~\cite{danz09,nabi13,wils14a,raja25}, for which we monitor the average voltage for the neurons, and if this crosses a threshold we input one cycle of the pre-computed input.  By design, this input will be a desynchronizing influence.  Another control input only occurs if the previous input has finished and the average voltage again crosses threshold, for example due to the synchronizing influence of the coupling.

\begin{figure}
    \centering
    \includegraphics[width=3.4in]{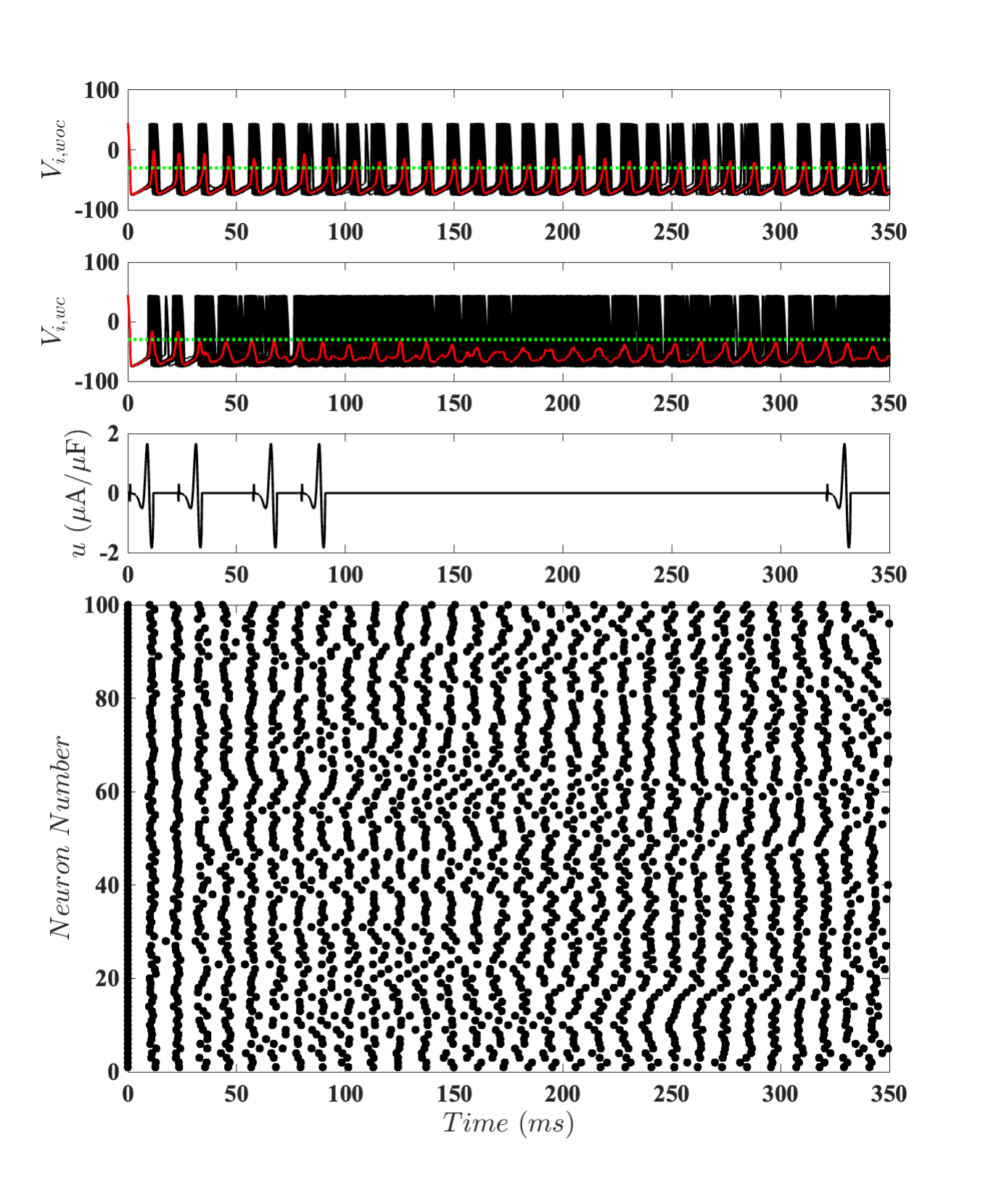}
    \caption{Results are shown for a population of $N = 100$ coupled neurons with $D = 2$ and homogeneous coupling strength $\alpha = 0.04$. 
The top panel shows the network's behavior without control, with synchronized dynamics. The second panel illustrates the same network with event-based control using $\tilde{u}^*(t)$, where the red trace shows the mean voltage, and the horizontal green dotted line shows the control activation threshold ($\bar{V} = -30$ mV). The third panel shows the control input. The bottom panel is a raster plot of the spike times of each neuron (each row) from which we see that the neurons are spiking at different times as our control signal is applied.}
    \label{fig:true}
    \end{figure}

\begin{figure}
    \centering
    \includegraphics[width=3.4in]{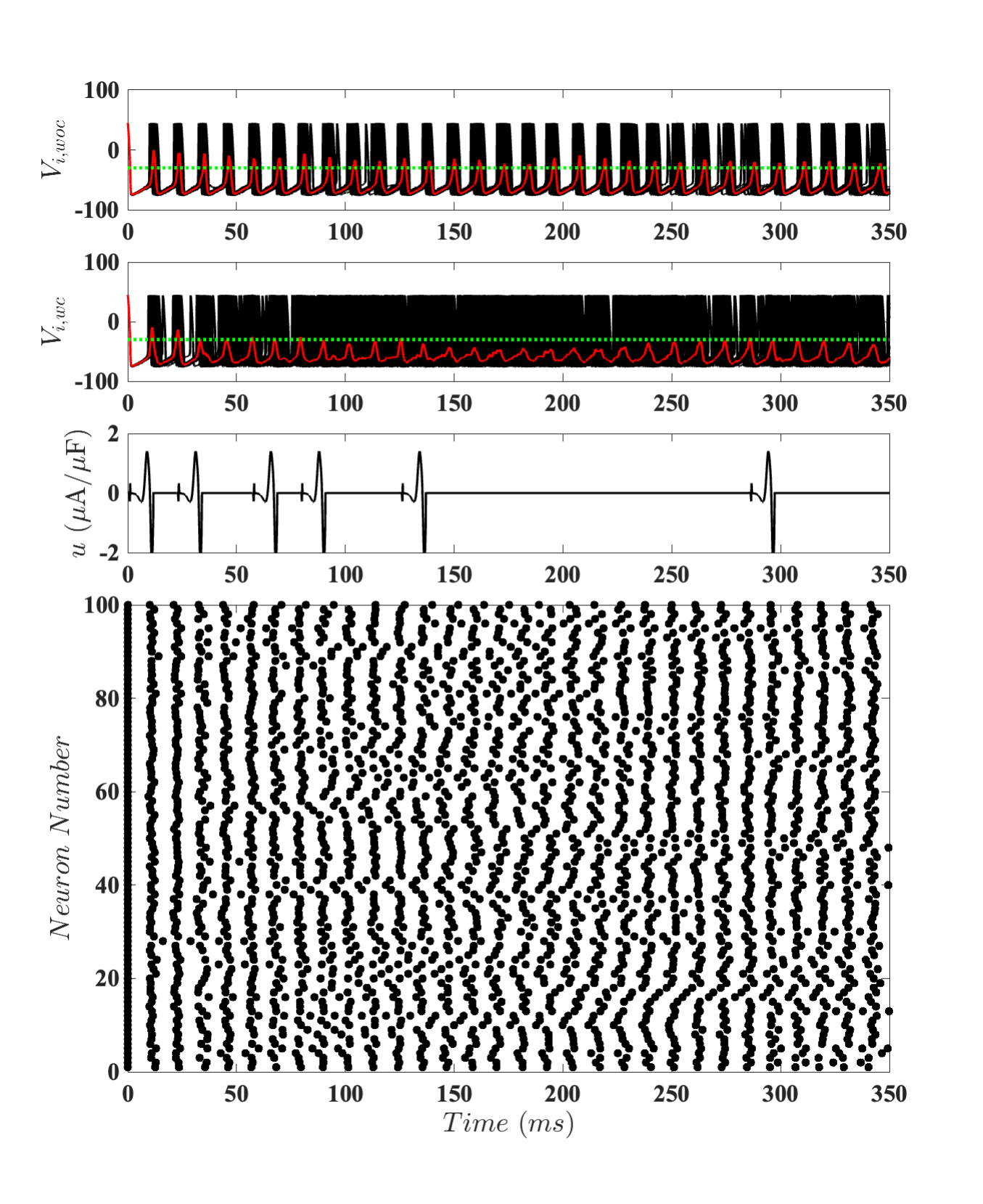}
    \caption{Results are shown for a population of $N = 100$ coupled neurons with $D = 2$ and homogeneous coupling strength $\alpha = 0.04$. 
The top panel shows the network's behavior without control, with synchronized dynamics. The second panel illustrates the same network with event-based control using $\tilde{u}^*_2(t)$, where the red trace shows the mean voltage, and the horizontal green dotted line shows the control activation threshold ($\bar{V} = -30$ mV). The third panel shows the control input. The bottom panel is a raster plot of the spike times, as for Figure~\ref{fig:true}
}. 
    \label{fig:Better}
    \end{figure}

Figures~\ref{fig:true} and~\ref{fig:Better} show results using the optimal input $u^*(t)$ and the nearly optimal input $u^*_2(t)$, respectively, as the input when the average voltage goes above threshold.  We note that these figures each show results for a particular realization of the noise terms $\eta_i(t)$.  Similar results are found for other noise realizations.  We see that both control schemes are effective: they can keep the neural population desynchronized.

% \begin{figure}
%    \centering
%    \includegraphics[width=3.2in]{Fig1.eps}
%    \caption{Results for Population-Level Simulations with $uopt_{true}$}
%    \label{fig:True}
%    \end{figure}

%     \begin{figure}
%    \centering
%    \includegraphics[width=3.2in]{Fig2.eps}
%    \caption{Results for Population-Level %Simulations with $uopt_{Zp}$}
%    \label{fig:Zp}
%    \end{figure}

One more cycle of input is required for Figure~\ref{fig:Better} compared with Figure~\ref{fig:true}.  This is true for the particular realizations of the noise for these figures, but won't always be the case.
In order to compare the necessary energy to keep the neurons desynchronized when using $u^*(t)$, $u^*_1(t)$, and $u^*_2(t)$ with this control scheme, we ran 100 noise realizations for each of these, and tracked the average energy used over the time interval shown in Figure~\ref{fig:Better}, namely from 0 to 350 ms.  Here energy is defined to be
%\[
$\int_0^{350} [u(t)]^2 dt$.
%\]
Results are shown in Table~\ref{table:2}.  Not surprisingly, the lowest average energy is found for the optimal input $u^*(t)$.   When $u^*_1(t)$ is used, this scheme requires approximately $26.5\%$ more energy.  On the other hand, when $u^*_2(t)$ is used, this scheme only requires approximately $5.8\%$ more energy.  This reinforces that $u^*_2(t)$ is a nearly optimal input for chaotic desynchronization.  Qualitatively similar results are found for other values of the noise and coupling strength.

\begin{table}[h!]
\centering
\caption{Mean and standard deviation of the energy from the average of 100 population-level simulations, with each simulation using a different realization of the noise.  }
\begin{tabular}{||c c c c||} 
 \hline
  & $u^*(t)$ & $u^*_1(t)$ & $u^*_2(t)$ \\ [0.5ex] 
 \hline\hline
  ${\rm mean} \left( \int u^2 \, dt \right)$ & 78.63 & 99.49 & 83.02 \\
 \hline
 ${\rm stdev} \left(\int u^2 \, dt \right)$ & 10.84 & 12.26 & 11.95 \\ [1ex] 
 \hline
\end{tabular}
\label{table:2}
\end{table}

\section{Synchronization Results} \label{sync}

In Section~\ref{desync} we compared the efficacy of the optimal input $u^*(t)$ and approximations $u^*_1(t)$ and $u^*_2(t)$ for desynchronizing pairs of neural oscillators.  By choosing a negative value of $\beta$, (\ref{G}) becomes
\begin{equation}
G[u(t)]=\int_0^{T} [u(t)^2 + |\beta| Z'(\theta(t))u(t)]dt.
\label{G_neg_beta}
\end{equation}
Minimizing (\ref{G}) then corresponds to simultaneously minimizing the energy and the Lyapunov exponent.  Typically, the minimal Lyapunov exponent will be negative, corresponding to synchronization.  Thus, by solving the two-point boundary value problem (\ref{thetadot}-\ref{findcontrol}) with negative $\beta$, we are finding the optimal input for {\it synchronization}.

\begin{figure}
    \centering
    \includegraphics[width=3.1in]{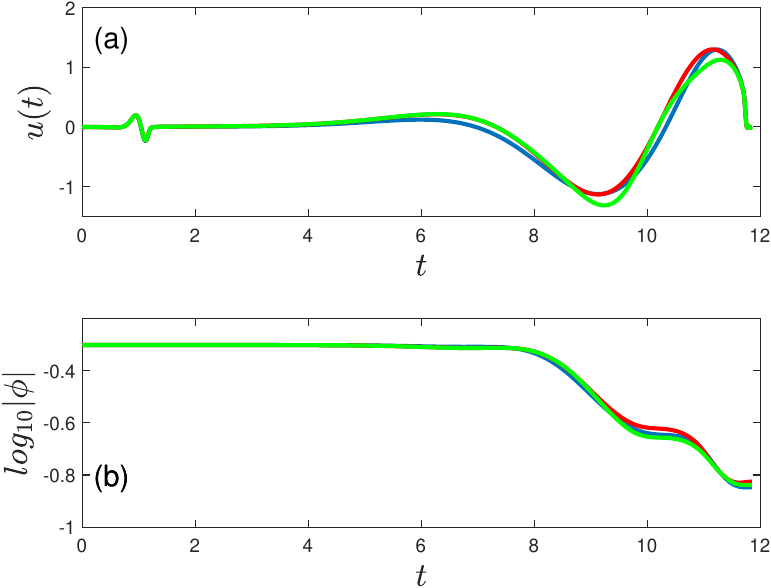}
    \caption{Results for the RHH PRC $Z_{\rm RHH}(\theta)$, with $\omega = 2 \pi/11.85$, and $\beta = -5$.  (a) shows the inputs $u^*(t)$ (blue), $u^*_1(t)$ (red), and $u^*_2(t)$ (green).  (b) shows the evolution of the phase difference $\phi$ for one cycle of each input, with colors as for (a).  The initial phase separation $\phi_0$ = 0.5.}
    \label{fig:RHH_sync}
    \end{figure}

Here we briefly explore our approximations for this optimal synchronization problem.  In particular, 
Figure~\ref{fig:RHH_sync} shows results for the PRC $Z_{\rm RHH}(\theta)$ for the RHH equations with $\beta = -5$, with (a) showing the inputs, and (b) showing the phase difference which results from these inputs with $\phi_0 = 0.5$.  We see that all three inputs $u^*(t), u^*_1(t)$, and $u^*_2(t)$ give very similar performance, with $u^*_2(t)$ giving a better approximation to the optimal performance than $u^*_1(t)$. (For $u^*(t), u^*_1(t)$, and $u^*_2(t)$, we get $\phi(T) = 0.142, 0.149$, and $0.145$, respectively.)

\section{Conclusion} \label{conclusion}

In this paper, we considered the desynchronization of neural oscillators by maximizing the Lyapunov exponent of the phase difference between pairs of oscillators, as motivated by deep brain stimulation treatment of neural disorders such as Parkinson's disease.  We derived two different approximations to optimal stimuli based on the oscillators' phase response curve.  Unlike the previous work in~\cite{wils14a}, these do not require numerical solution of a two-point boundary value problem; this might allow real-time control in which, for example, changes in $\beta$ are desired.  It was shown that one of these approximations, $u^*_2(t)$, has nearly optimal performance, and can be used with an event-based control scheme to desynchronize populations of noisy coupled neurons.

\section*{Appendix}

The neuron model which is considered in this paper is a two-dimensional reduction~\cite{{keen98,moeh06jmb}} of the Hodgkin-Huxley equations~\cite{hodg52d} given by (\ref{eq:Vndot}), with
\begin{eqnarray*}
f_V &=& \left[ I_b - \bar{g}_{Na}[m_{\infty}(V)]^3 (0.8 - n)(V - V_{Na}) \right. \\
&& - \left. \bar{g}_K n^4 (V - V_K) - \bar{g}_L (V - V_L) \right] / C, \label{eq:fv} \\
f_n &=& a_n(V)(1 - n) - b_n(V)n,
\label{eq:fn}
\end{eqnarray*}
\begin{eqnarray*}
m_{\infty}(V) &=& \frac{a_m(V)}{a_m(V)+b_m(V)} \;, \\
a_m(V) &=&  0.1(V+40)/(1-\exp(-(V+40)/10))  \; , \\
b_m(V) &=&  4\exp(-(V+65)/18)  \; , \\
a_n(V) &=&  0.01(V+55)/(1-\exp(-(V+55)/10))  \; , \\
b_n(V) &=&  0.125\exp(-(V+65)/80) \; , 
\end{eqnarray*}
\begin{eqnarray*}
&&V_{Na}=50 \; {\rm{mV}},\; V_K=-77 \; {\rm{mV}},\; V_L=-54.4 \; {\rm{mV}}, \\ 
&&\bar{g}_{Na}=120 \; {\rm{mS/cm^2}}, \; \bar{g}_K=36 \; {\rm{mS/cm^2}}, \\ 
&& \bar{g}_L=0.3 \; {\rm{mS/cm^2}}, \; C=1 \; {\rm{\mu F/cm^2}}.
\end{eqnarray*}

\noindent
Here \( I_b \) is the baseline current of the neuron.

\bibliographystyle{IEEEtran}
\bibliography{cdc25}

%\begin{thebibliography}{99}

%\bibitem{c1}
%J.G.F. Francis, The QR Transformation I, {\it Comput. J.}, vol. 4, 1961, pp 265-271.

%\bibitem{c2}
%H. Kwakernaak and R. Sivan, {\it Modern Signals and Systems}, Prentice Hall, Englewood Cliffs, NJ; 1991.

%\bibitem{c3}
%D. Boley and R. Maier, "A Parallel QR Algorithm for the Non-Symmetric Eigenvalue Algorithm", {\it in Third SIAM Conference on Applied Linear Algebra}, Madison, WI, 1988, pp. A20.

%\end{thebibliography}

\end{document}